\documentclass[12pt,leqno]{amsart}
\usepackage{amsfonts}

\newcommand{\Ff}{{\mathbb F}}
\newcommand{\Zz}{{\mathbb Z}}

\newcommand{\Qq}{{\mathbb Q}}

\newcommand{\Ll}{{\mathcal L}}
\newcommand{\Fl}{{\mathcal F}}

\newcommand\im{\operatorname{Im}}
\newcommand\Ker{\operatorname{Ker}}     %kernel
\newcommand\rank{\operatorname{rank}}   %rank
\newcommand\MM{{\overline M}}   %
\newcommand\pp{{\mathfrak p}}       %
\newcommand\bb{{\mathcal B}}        %
\newcommand\cc{{\mathcal C}}        %

\newtheorem{Theorem} {Theorem} [section]
\newtheorem{Proposition} [Theorem] {Proposition}
\newtheorem{Lemma} [Theorem] {Lemma}
\newtheorem{Corollary} [Theorem] {Corollary}

\newtheorem{Definition}[Theorem]{Definition}
\newtheorem{Remark}[Theorem]{Remark}

\newcommand{\Proof}{ \noindent{\bf Proof:}\quad }

\def\Im{\operatorname{Im}}
\def\PG{\operatorname{PG}}
\def\AG{\operatorname{AG}}

\def\QED{\qed\medskip\par}

\catcode`\@=11
\def\sideset#1#2#3{%
  \@mathmeasure\z@\displaystyle{#3}%
  \global\setbox\@ne\vbox to\ht\z@{}\dp\@ne\dp\z@
  \setbox\tw@\box\@ne
  \@mathmeasure4\displaystyle{\copy\tw@#1}%
  \@mathmeasure6\displaystyle{#3{#2}}%
  \dimen@-\wd6 \advance\dimen@\wd4 \advance\dimen@\wd\z@
  \hbox to\dimen@{}\mathop{\kern-\dimen@\box4\box6}%
}
\catcode`\@=12

\numberwithin{equation}{section}

\def\eqref#1{(\ref{#1})}

\begin{document}

\title[Invariant Factors of the Incidence Between Points and Subspaces]
{The Invariant Factors of the Incidence Matrices of points and subspaces in $\PG(n,q)$ and $\AG(n,q)$}

\author{David B. Chandler, Peter Sin, Qing Xiang}

\thanks{The second author was partially supported
by NSF grant DMS-0071060. The third author was partially supported by 
NSA grant MDA904-01-1-0036.}
\address{Department of Mathematical Sciences, University of Delaware,
  Newark, DE 19716, USA, email: chandler@math.udel.edu}
\address{Department of Mathematics, University of Florida, Gainesville, FL 32611, email: sin@math.ufl.edu}
\address{Department of Mathematical Sciences, University of Delaware,
  Newark, DE 19716, USA, email: xiang@math.udel.edu}

\keywords{}\subjclass{ 05E20 (Primary), 20G05 20C11 (Secondary)}
\renewcommand{\subjclassname}{\textup{2000} Mathematics Subject Classification}
\begin{abstract}
We determine the Smith normal forms of the incidence matrices 
of points and projective $(r-1)$-dimensional subspaces of $\PG(n,q)$
and of the incidence matrices of points and $r$-dimensional
affine subspaces of $\AG(n,q)$ for all $n$, $r$, and arbitrary prime power $q$.
\end{abstract}

\maketitle

\section{Introduction}
Let $\Ff_q$ be the finite field of order $q$, where $q=p^t$, $p$
is a prime and $t$ is a positive integer, and let $V$ be an
$(n+1)$-dimensional vector space over $\Ff_q$. We denote by
$\PG(V)$ (or $\PG(n,q)$ if we do not want to emphasize the
underlying vector space) the $n$-dimensional projective geometry
of $V$. The elements of $\PG(V)$ are subspaces of $V$ and two
subspaces are considered to be incident if one is contained in the
other. We call one-dimensional subspaces of $V$ \textit{points} of
$\PG(V)$ and we call $n$-dimensional subspaces of $V$
\textit{hyperplanes} of $\PG(V)$. More generally, we regard
$r$-dimensional subspaces of $V$ as projective $(r-1)$-dimensional
subspaces of $\mathrm{PG}(V)$. We will refer to $r$-dimensional
subspaces of $V$ as $r$-subspaces and denote the set of these
spaces in $V$ as $\mathcal L_r$. The set of projective points is
then $\mathcal L_1$. In this paper, we are concerned with the
incidence relation between $\mathcal L_r$ and $\mathcal L_1$.
Specifically, let $A$ be a (0,1)-matrix with rows indexed by
elements $Y$ of $\mathcal L_r$ and columns indexed by elements $Z$
of $\mathcal L_1$, and with the $(Y,Z)$ entry equal to 1 if and
only if $Z\subset Y$. We are interested in finding the Smith
normal form  (\cite{cohn}, p. 279) of $A$.

The incidence matrix $A$ has been studied at least since 1960s. In fact,
several authors have considered the more general incidence matrices $A_{r,s}$
of $r$-subspaces vs. $s$-subspaces, where $s$ is not necessarily one.
Most of their investigations were on the rank of $A_{r,s}$ over fields $K$ of
various characteristics. When $K=\Qq$, Kantor in \cite{bill} showed that
the matrix $A_{r,s}$ has full rank under certain natural conditions on $r$ and $s$, and when ${\rm char}(K)=\ell$, where
$\ell$ does not divide $q$, the rank of $A_{r,s}$ over $K$ was given by
Frumkin and Yakir \cite{fy}. The most interesting case is when
${\rm char}(K)=p$. In this case, the problem of finding the rank of $A_{r,s}$
is open in general (cf. \cite{godsil}). However, under the additional
condition $s=1$, Hamada \cite{ham} gave a complete solution to the problem of
finding the $p$-rank of $A$ (known as Hamada's formula). In this paper,
we are not only interested in the $p$-rank of $A$, but also the Smith normal form of
$A$ as an integral matrix. There are a couple of reasons for us to study this problem.

First, if we use the elements of $\mathcal  L_1$ as points and
use the elements of $\mathcal L_r$ as blocks, then we obtain what
is called a 2-design \cite{assmus} with ``classical parameters''.
It is known that there exist many 2-designs with
classical parameters \cite{bjl}. A standard way to distinguish
nonisomorphic designs with the same parameters is by comparing the
$p$-ranks of their incidence matrices. Unfortunately,
nonisomorphic designs sometimes have the same $p$-rank. In such
a situation, one can try to prove nonisomorphism of designs by
comparing the Smith normal forms of the incidence matrices
\cite{snf}. Therefore it is of interest to find Smith normal forms
of incidence matrices of designs.

Second, let $\Omega$ be an $n$-set. We say that an $r$-subset of $\Omega$
is {\it incident} with an $s$-subset of $\Omega$ if one is contained in the
other. In \cite{rick}, Wilson found a diagonal form of the
incidence matrix of $r$-subsets versus $s$-subsets of $\Omega$.
One can consider the $q$-analogue of this problem, namely, finding the
Smith normal form of the incidence matrix $A_{r,s}$ defined above.
So far we only succeeded in solving this problem in the case where $s=1$.
As far as we know, the problem of finding the $p$-rank of $A_{r,s}$ is open,
let alone finding the Smith form of $A_{r,s}$.

We briefly summarize previous work on or related to the problem of
finding the Smith form of the incidence between $\mathcal L_1$ and
$\mathcal L_r$. Hamada \cite{ham} determined the $p$-rank of the
incidence between projective points and $r$-subspaces of
$\PG(n,q)$ for any values of $p$, $t$, $r$, and $n$.  The work of
Hamada in \cite{ham} is based on an earlier paper of Smith
\cite{smith}.  A more conceptual proof of Hamada's formula, 
independent of \cite{smith}, was given in \cite{bsin}. Lander
\cite{lander} found the Smith form for the incidence between
points and lines in $\PG(2,q)$. Black and List \cite{black}
determined the invariant factors of the incidence between points
and hyperplanes in the case where $q=p$ (i.e., $t=1$). More
recently, Liebler \cite{liebler} and the second author 
each determined the
invariant factors of the incidence between points and hyperplanes
for general $q$. The invariant
factors of the incidence between points and arbitrary $r$-spaces
when $q=p$ (i.e., $t=1$) were computed in  \cite{sinp}. 
Finally, Liebler and the second author \cite{liebler}
had conjectured formulas for the invariant factors of the
incidence between points and arbitrary $r$-subspaces for general
$q$, and could prove their formulas in the cases where $q=p$,
$p^2$, or $p^3$. In this paper we use a combination of techniques
from number theory and representation theory to confirm this
conjecture.

In the following we will give a brief overview of the paper. For
convenience, we define the map
\begin{equation}\label{defeta}
\eta_{1,r}: \Zz^{\mathcal L_1}\rightarrow \Zz^{\mathcal L_r}
\end{equation}
by letting $\eta_{1,r}(Z)=\sum_{Y\in\mathcal L_r, Z\subset Y}Y$
for every $Z\in\mathcal L_1$, and then extending $\eta_{1,r}$
linearly to $\Zz^{\mathcal L_1}$. The matrix of $\eta_{1,r}$ with
respect to the basis $\mathcal L_1$ of $\Zz^{\Ll_1}$ and the basis
$\mathcal L_r$ of $\Zz^{\Ll_r}$ is exactly the matrix $A$ defined
above. We will use the same $\eta_{1,r}$ to denote the linear map
from $R^{\Ll_1}$ to $R^{\Ll_r}$ defined in the same way as above,
where $R$ is a certain $p$-adic local ring with maximal ideal $\pp$
and residue field $\Ff_q$ (see details in Section 2). The paper is
organized as follows. In Section 2, we introduce the monomial
basis ${\mathcal M}$ of $\Ff_q^{\Ll_1}$ and its Teichm\"uller
lifting to a basis ${\mathcal M}_R$ of $R^{\Ll_1}$. These bases
are very important for finding the Smith normal form of $A$. In
Section 3, we state our main theorem (Theorem~\ref{main}) which
gives the Smith normal form of $A$. We also include an elementary
proof of a well known fact stating that all the invariant factors
of $A$ are powers of $p$ except the last one. In Section 4, we
discuss Wan's theorem \cite{wan} on $p$-adic estimates of certain
multiplicative character sums. Wan's theorem can be applied directly
to our situation to give lower bounds on the ($p$-adic)
invariant factors of $A$ (now viewed as a matrix with entries from
$R$). In order to prove that these lower bounds indeed give the
$p$-adic invariant factors of $A$, considerable efforts are
needed. In Section 5, we prove that there exists a basis $\bb$ of
$R^{\Ll_1}$ whose reduction modulo $\pp$ is the monomial basis of
$\Ff_q^{\Ll_1}$ such that the matrix of $\eta_{1,r}$ with respect
to $\bb$ and some basis of $R^{\Ll_r}$ is the ($p$-adic) Smith
normal form of $A$. Next we prove a refinement of this result in
Section 6. We show that there exists a basis $\bb$ of $R^{\Ll_1}$
with the following properties:
\begin{enumerate}
\item $\bb$ contains certain elements of ${\mathcal M}_R$---we will
make this precise in Section 6;
\item the reduction modulo $\pp$ of $\bb$ is ${\mathcal M}$; and
\item there exists a basis $\cc$ of $R^{\Ll_r}$ such that the
matrix of $\eta_{1,r}$ with respect to $\bb$ and $\cc$ is the
$p$-adic Smith normal form of $A$.
\end{enumerate}
\noindent The proof of this result uses the natural action of the general
linear group on $R^{\Ll_1}$, Jacobi sums, and Stickelberger's
theorem on Gauss sums. In Section 7, combining the results
in previous sections, we give a proof of a more precise
statement (Theorem~\ref{pvalue}) which implies our main theorem. 
Finally in Section 8, we use our results in the projective geometry case 
to obtain the Smith normal form of the incidence matrix of points and 
$r$-flats of ${\rm AG}(n,q)$.

\section{Monomial Bases}

As we will see, most of the invariant factors of $A$ are $p$
powers. It will be helpful to view the entries of $A$ as coming
from some $p$-adic local ring. Let $q=p^t$ and let
$K=\Qq_p(\xi_{q-1})$ be the unique unramified extension of degree
$t$ over $\Qq_p$, the field of $p$-adic numbers, where $\xi_{q-1}$
is a primitive $(q-1)^{\rm th}$ root of unity in $K$. Let
$R=\Zz_p[\xi_{q-1}]$ be the ring of integers in $K$ and let
$\mathfrak p$ be the unique maximal ideal in $R$. Then $R$ is a
principal ideal domain, and the reduction of $R \
(\mathrm{mod}\,\mathfrak p)$ will be $\Ff_q$. Define $\bar x$ to
be $x\ (\mathrm{mod}\,\mathfrak p)$ for $x\in R$. Let $T_q$ be the
set of roots of $x^q=x$ in $R$ (a Teichm\"uller set) and let $T$
be the Teichm\"uller character of $\Ff_q$, so that $T(\bar x)=x$
for $x\in T_q$. We will use $T$ to lift a basis of
$\Ff_q^{\mathcal L_1}$ to a basis of $R^{\mathcal L_1}$.

In (\ref{defeta}), we defined the map $\eta_{1,r}$ from
$\Zz^{\Ll_1}$ to $\Zz^{\Ll_r}$. Now we use the same $\eta_{1,r}$
to denote the map from $R^{\mathcal L_1}$ to $R^{\mathcal L_r}$
sending a 1-space to the formal sum of all $r$-spaces incident
with it. The matrix $A$ is then the matrix of $\eta_{1,r}$ with
respect to the (standard) basis $\mathcal L_1$ of $R^{\mathcal
L_1}$ and the (standard) basis $\mathcal L_r$ of $R^{\mathcal
L_r}$. Crucial to our approach of finding the Smith form of $A$ is
what we call a monomial basis for $R^{\mathcal L_1}$. We introduce
this basis below.

We start with the monomial basis of $\Ff_q^{\mathcal L_1}$. This
basis was discussed in detail in \cite{bsin}. Let $V=\Ff_q^{n+1}$. Then
$V$ has a standard basis $v_0,v_1,\ldots ,v_{n}$, where
$$v_i=(\underbrace{0,0,\ldots ,0,1}_{i+1},0,\ldots ,0).$$
We regard $\Ff_q^V$ as the space of
functions from $V$ to $\Ff_q$. Any function $f\in \Ff_q^V$ can be
given as a polynomial function of $n+1$ variables corresponding to the $n+1$
coordinate positions: write the vector ${\mathbf x}\in V$ as
$${\mathbf x}=(x_0,x_1,\ldots ,x_n)=\sum_{i=0}^{n}x_iv_i;$$
then $f=f(x_0,x_1,\ldots ,x_n)$. The function $x_i$ is, for
example, the linear functional that projects a vector in $V$ onto
its $i^{\rm th}$ coordinate in the standard basis.

As a function on $V$, $x_i^q=x_i$, for each $i=0,1,\ldots ,n$, so we
obtain all the functions via the $q^{n+1}$ monomial functions
\begin{equation}\label{affbasis}
\{\prod_{i=0}^nx_i^{b_i}\mid 0\leq b_i<q, i=0,1,\ldots ,n\}.
\end{equation}
Since the characteristic function of $\{0\}$ in $V$ is
$\prod_{i=0}^n(1-x_i^{q-1})$, we obtain a basis for
$\Ff_q^{V\setminus\{0\}}$ by excluding $x_0^{q-1}x_1^{q-1}\cdots
x_n^{q-1}$ from the set in (\ref{affbasis}).

The functions on $V\setminus\{0\}$ which descend to $\mathcal L_1$
are exactly those which are invariant under scalar multiplication
by $\Ff_q^*$. Therefore we obtain a basis ${\mathcal M}$ of
$\Ff_q^{\mathcal L_1}$ as follows.
$${\mathcal M}=\{\prod_{i=0}^nx_i^{b_i}\mid 0\leq b_i<q,
\sum_i b_i\equiv 0\; ({\rm mod}\; q-1), (b_0,b_1,\ldots ,b_n)\neq
(q-1,q-1,\ldots ,q-1)\}.$$ This basis ${\mathcal M}$ will be
called the {\it monomial basis} of $\Ff_q^{\mathcal L_1}$, and its
elements are called {\it basis monomials}.

Now we lift the function $x_i: V\rightarrow \Ff_q$ to a function
$T(x_i): V\rightarrow R$, where $T$ is the Teichm\"uller character
of $\Ff_q$. For $(a_0,a_1,\ldots ,a_n)\in V$, we have
$$T(x_i)(a_0,a_1,\ldots ,a_n)=T(a_i)\in R.$$
For each basis monomial $\prod_{i=0}^nx_i^{b_i}$, we define
$T(\prod_{i=0}^nx_i^{b_i})$ similarly. We have the following
lemma.

\begin{Lemma}\label{nakayama}
The elements in the set
\begin{eqnarray*}
{\mathcal M}_R
&=&\big\{T(\prod_{i=0}^nx_i^{b_i})\mid 0\leq b_i<q,
 \sum_i b_i\equiv
0\; ({\rm mod}\; q-1),\\&& (b_0,b_1,\ldots ,b_n)\neq (q-1,q-1,\ldots
,q-1)\big\}
\end{eqnarray*}
 form a basis of the free $R$-module $R^{\mathcal L_1}$.
\end{Lemma}

\Proof To simplify notation, we use $M$ to denote the free
$R$-module $R^{\mathcal L_1}$, set $v=|{\mathcal L_1}|$, and
enumerate the elements of ${\mathcal M}_R$ as $f_1,f_2,\ldots
,f_{v}$. Since the images of the elements of ${\mathcal M}_R$ in
the quotient $M/{\mathfrak p}M$ are exactly the elements in
${\mathcal M}$, which form a basis of $\Ff_q^{\mathcal L_1}\cong
M/{\mathfrak p}M$ (as vector spaces over $\Ff_q$) , by Nakayama's
lemma \cite{am}, the elements in ${\mathcal M}_R$ generate $M$
and since their number equals $\rank M$, they form a basis.
\QED

The basis ${\mathcal M}_R$ will be called the {\it monomial
basis} of $R^{\mathcal L_1}$, and its elements are called {\it
basis monomials}.

\section{The Main Theorem}

Let $q=p^t$, and let $V$ be an $(n+1)$-dimensional space over
$\Ff_q$. As before we use $A$ to denote the $|\Ll_r|\times
|\Ll_1|$ matrix of the linear map $\eta_{1,r}:
\Zz^{\Ll_1}\rightarrow \Zz^{\Ll_r}$ with respect to the standard
bases of $\Zz^{\Ll_1}$ and $\Zz^{\Ll_r}$. It is known that all
invariant factors of $A$ (as a matrix over $\Zz$) are $p$-powers
except the last one, which is also divisible by $(q^r-1)/(q-1)$.
In \cite{sinp}, a proof was given using the structure of
the permutation module for ${\rm GL}(n+1,q)$ acting on
$\Ll_1$ over fields of characteristic prime to $p$. 
We give an elementary proof of the result.

\begin{Theorem}\label{p'part}
Let $A$ be the matrix of the map $\eta_{1,r}$ with respect to the standard
bases of $\Zz^{\Ll_r}$ and $\Zz^{\Ll_1}$, and let $v=|\Ll_1|$. The invariant
factors of $A$ are all $p$-powers except for the $v^\mathrm{th}$ invariant,
which is a $p$-power times $(q^r-1)/(q-1)$.
\end{Theorem}

\Proof We first define $\eta_{r,1}: \Zz^{\mathcal
L_r}\rightarrow\Zz^{\mathcal L_1}$ to be the linear map sending 
each element of $\Ll_r$ to the formal sum of all the 1-spaces incident with 
it. Then the matrix of $\eta_{r,1}$ with respect to the standard
bases of $\Zz^{\Ll_r}$ and $\Zz^{\Ll_1}$ is $A^{\top}$. For the purpose of proving this theorem, it will be more convenient to work with $A^{\top}$.

We define $$\epsilon:\Zz^{\mathcal L_1}\rightarrow\Zz$$ to be the 
map sending each element in $\Ll_1$ to 1. Clearly $\epsilon$ maps
$\Zz^{\mathcal L_1}$ onto $\Zz$ and $\im\eta_{r,1}$
onto $\left(\begin{smallmatrix}\frac{q^r-1}{q-1}\end{smallmatrix}\right)\Zz$. 
Thus, 
$$
\Zz^{\mathcal L_1}/(\Ker\epsilon+\im\eta_{r,1})\cong 
\Zz/\left(\begin{smallmatrix}\frac{q^r-1}{q-1}\end{smallmatrix}\right)\Zz.
$$
To finish the proof of the theorem, we are reduced to
showing that
$(\Ker\epsilon+\im\eta_{r,1})/\im\eta_{r,1}$ is a $p$-group.
We show that if  $x\in\Ker\epsilon$ 
then $q^{r-1}x\in\Im\eta_{r,1}$. Now $\Ker\epsilon$
is spanned by vectors of the form $u-w$, where $u$ and $w$ are vectors
representing individual elements in $\Ll_1$, so it is enough to
show that $q^{r-1}(u-w)$ is in $\im\eta_{r,1}$. Let $U$
be some $(r+1)$-subspace of $V$ which contains both $u$ and $w$.
We define $\tilde\eta_{1,r}$ to be the linear map which maps a
projective point to the formal sum of the $r$-subspaces which both
contain the point and are contained in $U$ and define $\mathbf
j_U$ to be the formal sum of all the projective points inside $U$.
Then $\eta_{r,1}$ restricted to $r$-subspaces inside $U$ and
$\tilde\eta_{1,r}$ are simply the hyperplane-to-point and
point-to-hyperplane maps for the space $U$. By standard
formula from design theory we have
$$\eta_{r,1}(\tilde\eta_{1,r}(z))=q^{r-1}z+\frac{q^{r-1}-1}{q-1}{\mathbf j}_U$$
for every $z\in \Ll_1$. Hence by setting $z=u$ and $z=w$
respectively, and subtracting the resulting equations, we get
$$\eta_{r,1}(\tilde\eta_{1,r}(u-w))=q^{r-1}(u-w)$$
which is the desired result.  \QED

In view of Theorem~\ref{p'part}, in order to get the Smith normal
form of $A$, we just need to view $A$ as a matrix with entries
from $\Zz_p$, the ring of $p$-adic integers, and get its Smith
normal form over $\Zz_p$. This will be the approach we take in the
rest of the paper. To state our main theorem, we need more
notation.

Let $\mathcal{H}$ denote the set of $t$-tuples
$\xi=(s_0,s_1,\ldots,s_{t-1})$ of integers satisfying (for $0 \le j
\le t-1$) the following:
\begin{equation}\label{H}\begin{array}{l}
(1)\quad 1 \le s_j \le n,
\\
(2)\quad 0 \le ps_{j+1}-s_j \le (p-1)(n+1),
\end{array}\end{equation} with the subscripts read (mod
$t$). The set ${\mathcal H}$ was introduced in \cite{ham}, and
used in \cite{bsin} to describe the module
structure of $\Ff_q^{\Ll_1}$ under the natural action of ${\rm
GL}(n+1,q)$. 

For a nonconstant basis monomial
$$f(x_0,x_1,\ldots,x_n)=x_0^{b_0}\cdots x_n^{b_n},$$
in ${\mathcal M}$, we expand the exponents
$$b_i=a_{i,0}+pa_{i,1}+\cdots+p^{t-1}a_{i,t-1}\quad 0\le
a_{i,j}\le p-1$$ and let
\begin{equation}\label{lambda}
\lambda_j=a_{0,j}+\cdots+a_{n,j}.
\end{equation}
Because the total degree $\sum_{i=0}^nb_i$ is divisible by $q-1$,
there is a uniquely defined $t$-tuple $(s_0,\ldots,s_{t-1})\in
\mathcal H$ \cite{bsin} such that
$$\lambda_j=ps_{j+1}-s_j.$$ Explicitly
\begin{equation}\label{defs}
s_j=\frac1{q-1}\sum_{i=0}^n\big(\sum_{\ell=0}^{j-1}p^{\ell+t-j}a_{i,\ell}+
\sum_{\ell=j}^{t-1}p^{\ell-j}a_{i,\ell} \big)
\end{equation}
One way of interpreting the numbers $s_j$ is that the total degree
of $f^{p^i}$ is $s_{t-i}(q-1)$, when the exponent of each
coordinate $x_i$ is reduced to be no more than $q-1$ by the
substitution $x_i^q=x_i$.  We will say that $f$ is of {\it type}
$\xi=(s_0,s_1,\ldots,s_{t-1})$. Also we say that the corresponding basis
monomial $T(f)\in {\mathcal M}_R$ is of {\it type} $\xi$.
(Note that in \cite{bsin} $\xi$ is called
{\it a tuple in $\mathcal H$} and the term {\it type}
is used for certain other $t$-tuples in bijection
with $\mathcal H$. However, since we will not use the latter there
is no risk of confusion.) 

Let $d_i$ be the coefficient of $x^{i}$ in the expansion of
$(\sum_{k=0}^{p-1}x^k)^{n+1}$. Explicitly,
$$d_i=\sum_{j=0}^{\lfloor i/p\rfloor}(-1)^j{n+1 \choose j}{n+i-jp \choose n}.$$

\begin{Lemma}
Let $d_i$ and $\lambda_j$ be as defined above. The number of basis
monomials in both ${\mathcal M}$ and ${\mathcal M}_R$ of type
$\xi=(s_0,s_1,\ldots,s_{t-1})$ is $\prod_{j=0}^{t-1}d_{\lambda_j}$.
\end{Lemma}
\Proof From (\ref{lambda}) each $\lambda_j$ is the sum of $n+1$
integers which can be anywhere from 0 to $p-1$.  The number of
such choices is the same as the coefficient of $x^{\lambda_j}$ in
$(\sum_{k=0}^{p-1}x^k)^{n+1}$. Counting the choices for each
$\lambda_j$ as $j$ runs from 0 to $t-1$ we get
$\prod_{j=0}^{t-1}d_{\lambda_j}$. \QED

We can now state the main theorem.

\begin{Theorem}\label{main}

Let $\mathcal{L}_1$ be the set of projective points and let
$\mathcal{L}_r$ be the set of projective $(r-1)$-spaces in
$\mathrm{PG}(n,q)$, and let $d_i$ and $\mathcal{H}$ be as above.
For each $t$-tuple $\xi=(s_0,s_1,\ldots,s_{t-1})\in \mathcal H$ let
$$\lambda_i=ps_{i+1}-s_i$$ and let
$$d_\xi=\prod_{i=0}^{t-1} d_{\lambda_i}.$$ Then
the $p$-adic invariant factors of the incidence matrix $A$ between
$\mathcal{L}_1$ and $\mathcal{L}_r$ are $p^\alpha$,
$0\leq \alpha\leq (r-1)t$, with multiplicity
$$m_\alpha=\sum_{\xi\in\mathcal H_\alpha}
d_\xi+ \delta(0,\alpha)$$ where \begin{equation}\label{halpha}
\mathcal H_\alpha = \Big\{(s_0,s_1,\ldots ,s_{t-1}) \in
\mathcal{H} \mid \sum_{i=0}^{t-1}\max\{0,r-s_i\}=\alpha\Big\},
\end{equation}
and
\begin{equation}
\delta(0, \alpha)=
\left\{
\begin{array}{ll}
1, & \mbox{if} \; \alpha=0,\\
0, & \mbox{otherwise}. \end{array}
\right.
\end{equation}
\end{Theorem}

\begin{Remark} {\rm (1). The theorem was conjectured by Liebler and the second author \cite{liebler}.

(2). The multiplicity of 1 among the $p$-adic invariant factors, $m_0$, is exactly the $p$-rank of $A$. From Theorem~\ref{main}, we have
$$m_0=1+\sum_{(s_0,s_1,\ldots ,s_{t-1})\in {\mathcal H}, s_i\geq r, \forall {i}}d_{(s_0,s_1,\ldots ,s_{t-1})}.$$
We mention that $d_{(s_0,\ldots,s_{t-1})}=d_{(n+1-s_0,\ldots,n+1-s_{t-1})}$ for each $(s_0,\ldots,s_{t-1})\in\mathcal H$, and  $d_0=1$, $d_i=0$ if $i<0$, and in fact $d_{(0,0, \ldots,0)}=d_{(n+1,n+1,\ldots,n+1)}=1$ and $d_\xi=0$ for all other cases that $\xi\not\in\mathcal H$. So the above $p$-rank formula is the same as the formula of Hamada \cite{ham}.

(3). We also mention that the largest $\alpha$ of the exponents of the $p$-adic invariant factors of $A$ is $(r-1)t$. It arises in the case where $\xi=(1,1,\ldots,1)$.  From Theorem~\ref{main}, we find that the multiplicity of $p^{(r-1)t}$ is
$$m_{(r-1)t}=d_{(1,1,\ldots ,1)}={n+p-1\choose n}^t,$$
which is one less than the $p$-rank of $\eta_{1,n}$.}
\end{Remark}

We indicate how we proceed to prove Theorem~\ref{main}. In
order to get the Smith normal form of $A$ over $R$, we will find two
invertible matrices $P$ and $Q^{-1}$ with entries in $R$, such that
$$A=PDQ^{-1},$$
where $D$ is a $|{\mathcal L}_r|\times |{\mathcal L}_1|$ diagonal
matrix with $p$ powers on its diagonal. The matrices
$Q$ and $P$ will come from basis changes in $R^{{\mathcal L}_1}$
and $R^{{\mathcal L}_r}$ respectively.

\vspace{0.1in}

Let $\{e_1,e_2,\ldots ,e_v\}$, where $v=|{\mathcal L}_1|$, be the
standard basis of $R^{{\mathcal L}_1}$, and let ${\mathcal
M}_R=\{f_1,f_2,\ldots ,f_{v}\}$ be the monomial basis of
$R^{{\mathcal L}_1}$ constructed in Lemma~\ref{nakayama}. For
$1\leq j\leq v$, let $f_j=\sum_{i=1}^v q_{ij}e_i$, $q_{ij}\in R$,
and let $Q=(q_{ij})$. Then
$$\eta_{1,r}(f_j)=\sum_{i=1}^v q_{ij}\eta_{1,r}(e_i).$$
Therefore the columns of $AQ$ are the vectors $\eta_{1,r}(f_j)$,
written with respect to the standard basis of $R^{{\mathcal
L}_r}$. For $1\leq j\leq v$, let $p^{a_j}$ be the largest power of
$p$ dividing every coordinate of $\eta_{1,r}(f_j)$. Then we try to
factorize $AQ$ as $PD$, where
$$D=
\begin{pmatrix}
p^{a_1} & 0 & 0 & \cdots & 0\\
0 & p^{a_2} & 0\\
0 & & \ddots & & \vdots\\
\vdots & & & p^{a_{v-1}}& 0\\
0 & & \cdots & 0 & p^{a_v}\\
0 & & \cdots & & 0\\
\vdots & & \ddots& & \vdots\\
0 & & \cdots & & 0
\end{pmatrix},$$
and $P$ is an $|\Ll_r|\times |\Ll_r|$ matrix whose first $v$
columns are ${\frac{1}{p^{a_j}}} \eta_{1,r}(f_j)$, $j=1,2,\ldots
,v$. In order to get the Smith normal form of $A$, we need to have
some information on $a_j$. For this purpose we need to have some
lower bound on the $p$-adic valuations of the coordinates of
$\eta_{1,r}(f_j)$. Let $f_j$ be a typical basis monomial
$T(x_0^{b_0}x_1^{b_1}\cdots x_n^{b_n})$ in ${\mathcal M}_R$, and
let $Y\in {\Ll_r}$. Then the $Y$-coordinate of $\eta_{1,r}(f_j)$
is
\begin{eqnarray*}
\eta_{1,r}(f_j)(Y)&=&\sum_{Z\subset Y, Z\in \Ll_1}f_j(Z)\\
                  &=&\frac 1{q-1}\sum_{{\mathbf x}\in \Ff_q^{n+1}\setminus\{(0,0,\ldots ,0)\}, {\mathbf x}\in Y}T^{b_0}(x_0)T^{b_1}(x_1)\cdots T^{b_n}(x_n),\\
\end{eqnarray*}
where in the last summation, ${\mathbf x}=(x_0,x_1,\ldots ,x_n)\in
\Ff_q^{n+1}$. Therefore the coordinates of $\eta_{1,r}(f_j)$ are
all multiplicative character sums. Thanks to a theorem of Wan
\cite{wan}, one can indeed obtain lower bounds on the $p$-adic
valuations of these multiplicative character sums. We discuss
Wan's theorem and its applications in the next section.

\section{Wan's Theorem}

We adopt the same notation as in Section 2. That is, $q=p^t$,
$K=\Qq_p(\xi_{q-1})$ is the unique unramified extension of degree
$t$ over $\Qq_p$, $R=\Zz_p[\xi_{q-1}]$ is the ring of integers in
$K$,  and $\mathfrak p$ is the unique maximal ideal in $R$. Define
$\bar x$ to be $x\ (\mathrm{mod}\,\mathfrak p)$ for $x\in R$. Let
$T_q$ be the set of roots of $x^q=x$ in $R$ and let $T$ be the
Teichm\"uller character of $\Ff_q$, so that $T(\bar x)=x$ for
$x\in T_q$. Then $T$ is a $p$-adic multiplicative character of
$\Ff_q$ of order $(q-1)$ and all multiplicative characters of
$\Ff_q$ are powers of $T$.  Following the convention of Ax
\cite{ax}, $T^0$ is the character that maps all elements of
$\Ff_q$ to 1, while $T^{q-1}$ maps 0 to 0 and all other elements
to 1.

For $0\le i\le n$ let $F_i(x_1,\ldots,x_r)$ be polynomials of degree $d_i$
over $\Ff_q$ and let $$\chi_i=T^{b_i}\quad (0\le b_i\le q-1)$$
 be multiplicative characters.  We want the $p$-adic valuation
$\nu_p(S_q(\chi,F))$ of the multiplicative character sum
$$S_q(\chi,F)=\sum_{\mathbf
x\in\Ff_q^r}\chi_0(F_0(\mathbf x))\cdots\chi_n(F_n(\mathbf x)).$$
For an integer $k\ge0$ we define $\sigma_q(k)$ to be the sum of
the digits in the expansion of $k$ as a base $q$ number and
$\sigma(k)$ as the sum of the digits in the expansion of $k$ as
a base $p$ number.  Wan's Theorem  (\cite{wan}, Theorem 3.1) is the
following:
\begin{Theorem} {\em (Wan).} Let $d=\max_i\,d_i$ and $q=p^t$.  Then the
$p$-adic valuation of $S_q(\chi,F)$ is at least
$$\sum_{\ell=0}^{t-1}\bigg\lceil\frac{r-\frac1{q-1} \sum_{i=0}^n
\sigma_q(p^{\ell}b_i)d_i}d\bigg\rceil.$$
\end{Theorem}

Here we state a slightly stronger version of the theorem, which
follows immediately from the proof in \cite{wan}:
\begin{Theorem}
$$\nu_p(S_q(\chi, F) )
\ge\sum_{\ell=0}^{t-1}\max\Big\{0,\bigg\lceil\frac{r-\frac1{q-1}
\sum_{i=0}^n
\sigma_q(p^{\ell}b_i)d_i}d\bigg\rceil\Big\}.$$
\end{Theorem}

We will use this theorem only in the case where each
$F_i$ is a linear homogeneous function. For the convenience of the reader
we specialize the proof given in \cite{wan}.

\begin{Theorem}\label{wanlow} For each $i$, $0\le i\le n$,
let $\overline F_i(\bar {\mathbf x})=\bar
\gamma_{i1}\bar x_1+\cdots+\bar \gamma_{ir}\bar x_r$ be a linear functional
on $\Ff_q^r$. Then
$$\nu_p(S_q(\chi,\overline F))\ge\sum_{\ell=0}^{t-1}\max\{0,r-\frac1{q-1}
\sum_{i=0}^n\sigma_q(p^{\ell}b_i)\}.$$
\end{Theorem}

\Proof
We will write
$$F_i({\mathbf x})=\gamma_{i1}\,x_1+\cdots+\gamma_{ir}\,x_r$$ to represent the
lifted functions from $T_q^r$ to $R$ with $\gamma_{ij}=T(\bar \gamma_{ij})$.
  Using the congruence $$T(\bar x)
\equiv x^{q^r}\ (\mathrm{mod}\ q^r)\ $$ for all $x\in R$ we get
\begin{equation}\label{sum}S_q(\chi,\overline F)\equiv
\sum_{{\bf x}\in
{T_q}^r}\big(F_0({\bf x})\big)^{b_0q^r}\cdots\big(F_n({\bf x})\big)^{b_nq^r}
\quad (\mathrm{mod}\ q^r).
\end{equation}
Expanding (\ref{sum}) we get
\begin{eqnarray*}
S_q(\chi,\overline F)&\equiv&
\end{eqnarray*}
\begin{equation}\label{expand}
 \sum_{
\begin{array}{c} k_{i1}+\cdots+k_{ir}=b_iq^r\\
0\le i\le n \end{array}}\prod_{i=0}^n
{{b_iq^r}\choose{k_{i1},
\ldots,k_{ir}}} \big(\prod_{i=0}^n\prod_{j=1}^r
{\gamma_{ij}}^{k_{ij}}\big)\big(
\prod_{j=1}^r\sum_{x\in T_q}x^{\sum_i k_{ij}}\big)
\quad(\mathrm{mod}\  q^r)
\end{equation}
We use the formula of Legendre, $\nu_p(k!)=(k-\sigma(k))/(p-1)$
and get that the $p$-adic valuation of the multinomial
coefficient part of (\ref{expand}) is
\begin{equation}\label{multin}
\frac1{p-1}\sum_{i=0}^n\big(b_iq^r-\sigma(b_i)-\sum_{j=1}^r
(k_{ij}-\sigma(k_{ij}))\big)=\frac1{p-1}\sum_{i=0}^n
\bigg(\sum_{j=1}^r\sigma(k_{ij})-\sigma(b_i)\bigg).
\end{equation}
For the Teichm\"uller set $T_q$ we have
\begin{equation}\label{tsum}
\sum_{x\in T_q} x^k=\left\{\begin{array}{ll} 0, &\mathrm{if}\ (q-1)\
\mathrm{does\ not\ divide}\ k,\\ q, & \mathrm{if}\  k=0, \\
q-1, & \mathrm{if}\ (q-1) |k\ \mathrm{and}\ k>0. \end{array}\right.
\end{equation}
Therefore, in (\ref{expand}) we only need to consider those terms
for which
\begin{equation}\label{ex1}
\sum_{i=0}^n k_{ij}\equiv 0\quad (\mathrm{mod}\, q-1)
\end{equation} for
all $j=1,2,\ldots ,r$.  Since $k\equiv \sigma_q(k)\ (\mathrm{mod}\ q-1),$ we
also have
\begin{equation}\label{ex2}
\sum_{i=0}^n\sigma_q(k_{ij})\equiv 0\quad
(\mathrm{mod}\, q-1).
\end{equation}
Given $k_{ij}$ such that $\sum_{j=1}^rk_{ij}=b_iq^r$ for $0\leq i\leq n$ and (\ref{ex1}) is satisfied, assume that $s$ coordinates of the vector
$$(\sum_{i=0}^nk_{i1},\sum_{i=0}^nk_{i2},\ldots ,\sum_{i=0}^nk_{ir})$$
are not identically 0. Then the same is true for
the corresponding entries of the vector
\begin{equation}\label{vec}
(\sum_{i=0}^n\sigma_q(k_{i1}),\sum_{i=0}^n\sigma_q(k_{i2}),\ldots ,\sum_{i=0}^n\sigma_q(k_{ir})).
\end{equation}
Summing up the entries of the vector in (\ref{vec}) we get
\begin{equation}\label{inequal}
s(q-1)-\sum_{i=0}^nb_i\le\sum_{i=0}^n\bigg(\sum_{j=1}^r
\sigma_q(k_{ij})-b_i\bigg).
\end{equation}
We note that for a non-negative integer $\ell$,
(\ref{ex2}) still holds with $\sigma_q(k_{ij})$ replaced by
$\sigma_q(p^{\ell}k_{ij})$.
Also $\sum_{i=0}^n\sigma_q(p^{\ell}k_{ij})$ is not identically 0 for the
same $s$ subscripts of $j$.  Thus we have
$$s(q-1)-\sum_{i=0}^n\sigma_q(p^{\ell}b_i)
\le\sum_{i=0}^n\bigg(\sum_{j=1}^r
\sigma_q(p^{\ell}k_{ij})-\sigma_q(p^{\ell}b_i)\bigg).$$
Noting that the right-hand side is non-negative since $\sum_{j=1}^rk_{ij}=b_iq^r$, we sum over $\ell$ to get
$$\sum_{\ell=0}^{t-1}\max
\big\{0,s(q-1)-\sum_{i=0}^n\sigma_q(p^{\ell}b_j)\big\}
\le\frac{q-1}{p-1}\sum_{i=0}^n\bigg(\sum_{j=1}^r
\sigma(k_{ij})-\sigma(b_i)\bigg),$$
using the fact that
$$\sum_{\ell=0}^{t-1}\sigma_q(p^{\ell}k)=\frac{q-1}{p-1}\sigma(k).$$
Comparing with  (\ref{multin}) we get that
each term of (\ref{expand}) (with $k_{ij}$ satisfying (\ref{ex1})) has
$p$-adic valuation at least
$$t(r-s)+\sum_{l=0}^{t-1}\max\bigg\{0,s-\frac1{q-1}
\sum_{i=0}^n\sigma_q(p^{\ell}b_i)\bigg\}\ge
\sum_{l=0}^{t-1}\max\bigg\{0,r-\frac1{q-1}
\sum_{i=0}^n\sigma_q(p^{\ell}b_i)\bigg\}.$$
This completes the proof. \QED

We now apply Wan's theorem to our situation. Let
$f=T(x_0^{b_0}x_1^{b_1}\cdots x_n^{b_n})\in {\mathcal M}_R$ be a basis monomial. We
use Theorem~\ref{wanlow} to give a lower bound on the $p$-adic valuation of the
coordinates of $\eta_{1,r}(f)$. Note that the coordinates of $\eta_{1,r}(f)$ are indexed by the
$r$-spaces in $\Ll_r$. An $r$-subspace $Y$ of $V=\Ff_q^{n+1}$ can
be defined by a system of $(n+1-r)$ independent linear homogeneous
equations. Putting the $n+1-r$ equations in reduced row echelon
form, we have $r$ coordinates which can run freely through $\Ff_q$
and the remaining $n+1-r$ coordinates are linear functions of
those $r$ coordinates.  Without loss of generality we label the
free coordinates $(x_0,\ldots,x_{r-1})=\mathbf x$ and express the
defining equations of $Y$ as $x_i=F_i(\mathbf x)$ for $(r\le i\le
n)$. The $Y$-coordinate of $\eta_{1,r}(f)$ is:
\begin{equation}\label{coordinate}
\eta_{1,r}(f)(Y)=\frac1{q-1}\sum_{\mathbf x\in {\Ff_q}^r\setminus\{(0,0,\ldots ,0)\}}T^{b_0}(x_0)\cdots T^{b_{r-1}}(x_{r-1})
T^{b_r}(F_r(\mathbf x))\cdots T^{b_n}(F_n(\mathbf x)).
\end{equation}

\begin{Lemma}\label{lower}
Let $f(x_0,\ldots,x_n)=T({x_0}^{b_0}\cdots{x_n}^{b_n})$ be a nonconstant
basis monomial in ${\mathcal M}_R$.  Then every coordinate
of the image vector $\eta_{1,r}(f)$ is divisible by $p^\alpha$ with
\begin{equation}\label{alpha}
\alpha=\sum_{i=0}^{t-1}\max\{0,r-s_i\},
\end{equation}
where $(s_0,s_1,\ldots ,s_{t-1})$ is the type of $f$ as defined
in (\ref{defs}).
\end{Lemma}
\Proof Let $Y$ be an arbitrary $r$-space in $\Ll_r$. By the above
discussion, we assume that $Y$ is defined by $x_i=F_i({\mathbf
x})$, $i=r,r+1,\ldots ,n$ and ${\mathbf x}=(x_0,x_1,\ldots
,x_{r-1})\in \Ff_q^r$. The $Y$-coordinate of $\eta_{1,r}(f)$ is
then given by (\ref{coordinate}), and by Theorem~\ref{wanlow}, we
have
$$\nu_p(\eta_{1,r}(f)(Y))\geq \sum_{\ell=0}^{t-1}\max\{0,r-\frac1{q-1}
\sum_{i=0}^n\sigma_q(p^{\ell}b_i)\}.$$ Recalling that the type of
$f$ is denoted by $(s_0,s_1,\ldots ,s_{t-1})$ and noting that
$s_{t-\ell}=\frac1{q-1}\sum_{i=0}^n \sigma_q(p^{\ell}b_i)$
(reading $s_t$ as $s_0$) for all $\ell=0,1,\ldots ,t-1$, we have
$$\nu_p(\eta_{1,r}(f)(Y))\ge\sum_{i=0}^{t-1}\max\{0,r-s_i\}.$$
This completes the proof. \QED

Let $Q$ be the basis change matrix between the standard basis and
the monomial basis ${\mathcal M}_R=\{f_1,f_2,\ldots ,f_v\}$ of
$R^{\Ll_1}$ (as used in Section 3). Using Lemma~\ref{lower}, we
see that one can factorize $AQ$ as $PD$, where
$$D=
\begin{pmatrix}
p^{\alpha_1} & 0 & 0 & \cdots & 0\\
0 & p^{\alpha_2} & 0\\
0 & & \ddots & & \vdots\\
\vdots & & & p^{\alpha_{v-1}}& 0\\
0 & & \cdots & 0 & p^{\alpha_v}\\
0 & & \cdots & & 0\\
\vdots & & \ddots& & \vdots\\
0 & & \cdots & & 0
\end{pmatrix},$$
$p^{\alpha_i}$ corresponds to the basis monomial $f_i\in {\mathcal
M}_R$ with type $(s_0,s_1,\ldots ,s_{t-1})$,
$$\alpha_i=\sum_{j=0}^{t-1}{\rm max}\{0,r-s_j\},$$
and $P$ is an $|\Ll_r|\times |\Ll_r|$ matrix whose first $v$
columns are ${\frac{1}{p^{\alpha_i}}} \eta_{1,r}(f_i)$,
$i=1,2,\ldots ,v$. We still need to show that $D$ (with the
diagonal entries suitably arranged) is indeed the Smith normal
form of $A$.

\section{$p$-Filtrations and Smith Normal Form Bases}

Let $R=\Zz_p[\xi_{q-1}]$ with maximal ideal $\pp=pR$ and residue
field $\Ff_q$, and let $\eta_{1,r}: R^{\mathcal L_1}\rightarrow
R^{\mathcal L_r}$ be the map defined before. In this section we
prove that there exists a basis $\bb$ of $R^{\Ll_1}$ whose
reduction modulo $\pp$ is the monomial basis of $\Ff_q^{\Ll_1}$
such that the matrix of $\eta_{1,r}$ with respect to $\bb$ and
some basis of $R^{\Ll_r}$ is the Smith normal form of
$\eta_{1,r}$. We begin with some general results on injective
homomorphisms of free $R$-modules. 

For any free $R$ module $M$ we set $\MM=M/\pp M$
and for any $R$-submodule $L$ of $M$, let $\overline
L=(L+\pp M)/\pp M$ be the image in $\MM$.

Let $\phi:M\rightarrow N$ be an injective homomorphism of free
$R$-modules of finite rank, with $\rank M=m\geq 1$.

Let
$$
N'=\{x\in N\mid \exists j\geq 0, p^jx\in\im\phi\}
\label{N'}
$$
Then $N'$ is the smallest $R$-module direct summand of 
$N$ containing $\im\phi$, (sometimes called its {\it purification})
and is also of rank $m$.  The invariant factors of $\phi$
stay the same if we change the codomain to $N'$. This 
will often allow us to reduce to the case $\rank N=m$.

Define
$$
M_i=\{m\in M \ \mid \ \phi(m)\in p^iN\}, \quad i=0,1,...
$$

Then we have a filtration
$$
M=M_0\supseteq M_1\supseteq \cdots
$$
of $M$
and the filtration
$$
\MM=\MM_0\supseteq \MM_1\supseteq \cdots
$$
of $\MM$.

Since $\phi$ is injective, and $N'/\im\phi$ has finite
exponent, it follows that there exists a
smallest index $\ell$ such that $\MM_{\ell}=0$. 
So we have a finite  filtration
$$
\MM=\MM_0\supseteq \MM_1\supseteq\cdots\supseteq \MM_{\ell}=\{0\}.
$$
Note that the inclusions need not be strict, though
the last one is, by minimality of $\ell$.

\begin{Proposition}
For $0\leq i\leq {\ell-1}$, $p^i$ is an invariant factor of $\phi$
with multiplicity $\dim (\MM_i/\MM_{i+1})$.
\end{Proposition}
\Proof 
The theory of modules over PIDs says there are bases of $M$
and $N'$ such that $\phi$ is represented by an $m\times m$
diagonal matrix whose entries are the invariant factors
of $\phi$. From this matrix we see that the multiplicity
of $p^i$  is $\dim (\MM_i/\MM_{i+1})$. \QED

Let us start with a basis ${\overline\bb}_{\ell-1}$ of
$\MM_{\ell-1}$ and extend it to a basis  of $\MM_{\ell-2}$ by
adding a set ${\overline\bb}_{\ell-2}$ of vectors and so on until
we have a basis
$$
{\overline\bb}={\overline\bb}_0\cup{\overline\bb}_1\cup\ldots\cup{\overline\bb}_{\ell-1}
$$
of $\MM$. At each stage we also select a set $\bb_i\subset M_i$ of
preimages of ${\overline\bb}_i$ and expand the sets in the same
way. The resulting set $\bb=\cup_{i=0}^{\ell-1}\bb_i$ is a basis
of $M$, by Nakayama's lemma. 

We show that this basis can be used
to compute the Smith normal form of $\phi$, namely that there is a
basis $\cc$ of $N$ such that the matrix of $\phi$ with respect to
$\bb$ and $\cc$ is the Smith normal form. 

Now for $e$ in $\bb_i$, we have $p^i\parallel\phi(e)$,
so $y=\frac1{p^i}\phi(e)$ is an element of $N'$.
The elements $y$ thus obtained
from all elements of $\bb$ are linearly independent elements
of $N'$, since $\phi$ is injective. Moroever, the index of 
$\im\phi$ in the $R$-submodule of $N'$ generated 
by these elements $y$ is equal to the index of $\im\phi$
in $N'$ by the proposition. Therefore, these elements
$y$ form a basis of $N'$. The matrix of
$\phi$ with respect to $\bb$ and any basis of $N$
obtained by extending this basis will then be in Smith 
normal form.

For convenience, we introduce a special name for bases
such as $\bb$ above.

\begin{Definition} {\rm We will call a basis $\bb$ of $M$
an {\it SNF basis of $M$ for $\phi$} if 
$\bb=\cup_{i=0}^{\ell -1}\bb_i$, where for each $i$
we have  $\bb_i\subseteq M_i$ and $\bb_i$ maps bijectively to a basis
of $\MM_i/\MM_{i+1}$ under the composite map $M_i\rightarrow \MM_i
\rightarrow \MM_i/\MM_{i+1}$.}
\label{SNFbasis}
\end{Definition}

We now apply the above general theory to our situation. We will
look at the case where $M=R^{\Ll_1}$, $N=R^{\Ll_r}$,
$\phi=\eta_{1,r}$. Let $G={\rm GL}(n+1,q)$. Then 
$G$ acts on $\Ll_1$ and $\Ll_r$ and the map
$\eta_{1,r}$ is an injective homomorphism of 
$RG$-modules, so  the $M_i$ are $RG$-modules
and the $\MM_i$ are $\Ff_qG$-modules.

We will use the following special properties 
of the $\Ff_qG$-module $\Ff_q^{\Ll_1}$.

\begin{Proposition}
\label{AB}
\begin{enumerate}
\item Two basis monomials of the same type generate the same
$\Ff_qG$-submodule of $\Ff_q^{\Ll_1}$.
\item Every $\Ff_qG$-submodule of $\Ff_q^{\Ll_1}$ has a basis
consisting of of all basis monomials in the submodule.
\end{enumerate}
\end{Proposition}
\begin{proof}
Part (1) is immediate from  \cite{bsin}, Theorem B. 
( The field in
\cite{bsin} is taken to be an algebraically closed field $k$,
 not $\Ff_q$, but it follows from \cite{bsin}, Theorem A
that in fact all the $kG$-submodules of $k^{\Ll_1}$ are
simply scalar extensions of $\Ff_qG$-submodules of
$\Ff_q^{\Ll_1}$, so for example \cite[Theorems A, B]{bsin} hold
also over $\Ff_q$.)
Let $S$ be an $\Ff_qG$-submodule of $\Ff_q^{\Ll_1}$ and
let ${\mathcal K}\subseteq {\mathcal H}\cup \{(0,\ldots,0)\}$
be the set of tuples of the composition factors of $S$.
Let $S'$ be the $\Ff_qG$-submodule generated by all 
basis monomials with tuples
in $\mathcal K$. By \cite{bsin}, Theorem B, $S'$ is the smallest
$\Ff_qG$-submodule such that the set of tuples
of its composition factors contains $\mathcal K$, so
$S'=S$. Hence, by \cite{bsin}, Theorem B,  in the expression of any
element of $S$ as a linear combination of basis monomials,
only basis monomials with tuples in $\mathcal K$ occur,
proving (2).
\end{proof}

\begin{Corollary}\label{snfbasis} $R^{\Ll_1}$ has an SNF basis
for $\eta_{1,r}$
whose image in $\Ff_q^{\Ll_1}$ is the monomial basis.
\end{Corollary}

\Proof By Proposition~\ref{AB}(2) we can choose ${\overline\bb}_i$
in the construction above to be the set of monomials in 
$\MM_i$ which are not in $\MM_{i+1}$.
\QED

Whenever we have a basis $\bb$ of $R^{\Ll_1}$ whose reduction modulo
$\pp$ is the monomial basis, the type of an element
of $\bb$ will always mean the type of its image in the
monomial basis.

\begin{Corollary}\label{types}
Let $\bb$ be an SNF basis of $R^{\Ll_1}$ for $\eta_{1,r}$
 whose image in
$\Ff_q^{\Ll _1}$ is the monomial basis. Then the invariants
corresponding to two elements of $\bb$ of the same type are equal.
\end{Corollary}

\Proof Let $e$, $f\in\bb$ be two such basis elements, with images
$\overline e$ and $\overline f$. Then
\begin{eqnarray*}
e\in M_j&\iff & \overline e\in \MM_j\quad \text{(def. of SNF basis)}\cr
&\iff & \overline f\in \MM_j\quad \text{(Proposition~\ref{AB}(1))}\cr 
&\iff & f\in M_j\quad \text{(def. of SNF basis)}.\cr
\end{eqnarray*}
\QED

\section{Jacobi Sums and The Action of The General Linear Group on $R^{\Ll_1}$}

In this section we will prove a refinement of Corollary~\ref{snfbasis} (see Lemma~\ref{basis} for 
details). In order to prove this refinement, we need to use Jacobi sums and 
the action of the general linear group on $R^{\Ll_1}$. We first define Jacobi sums.

Let $T$ be the Teichm\"uller character of $\Ff_q$ defined in Section 2,
where $q=p^t$. We know that $T$ is a $p$-adic multiplicative character of
$\Ff_q$ of order $(q-1)$ and all multiplicative characters of
$\Ff_q$ are powers of $T$.  Again we adopt the convention that
$T^0$ is the character that maps all elements of $\Ff_q$ to 1,
while $T^{q-1}$ maps 0 to 0 and all other elements to 1.

For any two integers $b_0$ and $b_1$, we define
\begin{equation}\label{defjacobi}
J(T^{b_0},T^{b_1})=\sum_{x_0\in \Ff_q}T^{b_0}(x_0)T^{b_1}(1-x_0)
\end{equation}
From the above definition and our convention on $T^0$ and
$T^{q-1}$, we see that if $b_0\not\equiv 0$ (mod $q-1$), then
$$J(T^{b_0},T^0)=0, \;{\rm and}\; J(T^{b_0},T^{q-1})=-1.$$
Also we have $J(T^{-1}, T)=1$. The Jacobi sum $J(T^{b_0},T^{b_1})$
lies in $R=\Zz_p[\xi_{q-1}]$. Naturally we want to know its
$p$-adic valuation. Using Stickelberger's
theorem on Gauss sums \cite{stick} (see \cite{dwork} for further
reference) and the well known relation between Gauss and Jacobi
sums, we have

\begin{Theorem}\label{jacobidiv}
Let $b_0$ and $b_1$ be integers such that $b_i\not\equiv 0$ (mod $q-1$), 
$i=0, 1$, and $b_0+b_1\not\equiv 0$ (mod $q-1$). For any integer $b$, 
we use $\sigma(b)$ to denote the sum of digits in the expansion of the 
least nonnegative residue of $b$ modulo $q-1$ as a base $p$ number. Then
$$\nu_p(J(T^{-b_0},T^{-b_1}))=\frac{\sigma(b_0)+\sigma(b_1)-\sigma(b_0+b_1)}{p-1}.$$
In other words, the number of times that $p$ divides $J(T^{-b_0},T^{-b_1})$ is equal to the number of carries in the addition $b_0+b_1$ (mod $q-1$).
\end{Theorem}

We will now construct an element of $RG$ with certain special properties. For this purpose, we will first describe 
the action of $G$ on $R^{\Ll_1}$. 
We think of elements of $\Ll_1$ in homogeneous
coordinates as row vectors and elements of $G$
as matrices acting by right multiplication.
Then $R^{\Ll_1}$ is the left $RG$-module given in the following way.
For each function $f\in R^{\Ll_1}$, the function $gf$ is given by
$$
(gf)(Z)=f(Zg),\qquad Z\in\Ll_1.
$$

Let $f_i=T(x_0^{b_0}x_1^{b_1}\cdots x_n^{b_n})\in \mathcal M_R$ be an
arbitrary basis monomial.
Let $\xi=\xi_{q-1}$, a primitive $(q-1)^{\rm th}$
root of unity in the Teichm\"uller set $T_q\subset R$, and 
let $\bar\xi$
its reduction modulo $p$. We define $g_\ell\in G$ to be the 
element which replaces $x_0$ by $x_0+\bar\xi^{\ell}x_1$
and leaves all other $x_i$ unchanged.
%$$g_\ell=
%\begin{pmatrix}
%1 & 0 & 0 & \cdots & 0\\
%\overline{t}^{\ell} & 1 & 0\\
%0 & & \ddots & & \vdots\\
%\vdots & & & 1& 0\\
%0 & & \cdots & 0 & 1\\
%\end{pmatrix}.$$
Then
$$g_\ell f_i=T\big((x_0+{\bar\xi}^\ell
x_1)^{b_0}x_1^{b_1}\cdots x_n^{b_n}\big).$$ Let
$g=\sum_{\ell=0}^{q-2}\xi^{-\ell}g_\ell\in RG$. The following lemma
gives us $gf_i$.

\begin{Lemma}\label{specialg}
Let $f_i$ and $g$ be as given.  Then
$$gf_i=\left\{\begin{array}{cl}
0, &\mathrm{if}\ b_0=0\\
T(x_0^{q-2}x_1^{b_1+1}x_2^{b_2}\cdots
x_n^{b_n}), &\mathrm{if}\ b_0=q-1\\
\big(q(1-T(x_0^{q-1}))-1\big)T(x_1^{b_1+1}x_2^{b_2}\cdots
x_n^{b_n}),&\mathrm{if}\ b_0=1\\
-J(T^{-1},T^{b_0})T(x_0^{b_0-1}x_1^{b_1+1}x_2^{b_2}\cdots
x_n^{b_n}),&\mathrm{otherwise.}
\end{array}\right.$$
\end{Lemma}

\Proof First note that
\begin{eqnarray*}
J(T^{-1},T^{0})&=&0,\\ J(T^{-1},T^{q-1})&=&-1,
\end{eqnarray*}
so the cases $b_0=0$ and $b_0=q-1$ are really covered by the
general case. Therefore we will only consider two cases.

\vspace{0.1in}
\noindent{\bf Case 1.}  $b_0\neq 1$. First assume that $x_0$ and
$x_1$ are both nonzero. We have
\begin{eqnarray}
gf_i&=&\sum_{\ell=0}^{q-2}\xi^{-\ell}g_\ell f_i\\ \label{eq82}
&=&T(x_1^{b_1}\cdots
x_n^{b_n})\sum_{\ell=0}^{q-2}T^{-1}(\bar\xi^{\ell})T^{b_0}(x_0+\bar\xi^\ell
x_1)\\ &=&
T(x_1^{b_1}\cdots
x_n^{b_n})\sum_{u \in\Ff_q}T^{-1}(-\frac{x_1u}{x_0})
T^{b_0}\bigg(1-(-\frac{x_1u}{x_0})
\bigg)T(-1)T(x_0^{b_0-1}x_1)\\ \label{eq84}&=&
-J(T^{-1},T^{b_0})T(x_0^{b_0-1}x_1^{b_1+1}x_2^{b_2}\cdots
x_n^{b_n}).
\end{eqnarray}
If $x_1=0$ we verify directly that (\ref{eq82}) and (\ref{eq84})
are both zero, so the formula is still valid. If $x_0=0$, since
$b_0\ne1$, we see that (\ref{eq84}) is 0; and (\ref{eq82}) is also
0, since a nontrivial (multiplicative) character summed over
$\Ff_q$ is zero. Therefore the formula still holds.

\vspace{0.1in}
\noindent{\bf Case 2.} $b_0=1$. In this case
$$gf_i=T(x_1^{b_1}x_2^{b_2}\cdots x_n^{b_n})\sum_{\ell=0}^{q-2}T(\bar\xi^{-\ell}x_0+x_1).$$
If $x_0=0$, then $gf_i=(q-1)T(x_1^{b_1+1}x_2^{b_2}\cdots
x_n^{b_n})$. If $x_0\neq 0$ but $x_1=0$, then clearly we have
$gf_i=0$. If $x_0\neq 0$ and $x_1\neq 0$, then using the same
calculations as in the case $b_0\neq 1$, we have
$$gf_i=-J(T^{-1},T)T(x_1^{b_1+1}x_2^{b_2}\cdots x_n^{b_n})=-T(x_1^{b_1+1}x_2^{b_2}\cdots x_n^{b_n}).$$
In summary, the formula for $gf_i$ in this case is
$$\big(q(1-T(x_0^{q-1}))-1\big)T(x_1^{b_1+1}x_2^{b_2}\cdots
x_n^{b_n}).$$
This completes the proof.\QED

\begin{Corollary}\label{frobenius}
Let $f_i=T(x_0^{b_0}x_1^{b_1}\cdots x_n^{b_n})\in \mathcal M_R$ be
a basis monomial, and let $g(j)=\sum_{\ell=0}^{q-2}\xi^{-\ell
}g_{\ell p^{-j}}$ be the $j^{\mathrm{th}}$ Frobenius analog of $g$
in Lemma~\ref{specialg} above. Then
$$g(j)f_i=\left\{\begin{array}{cl}
0, &\mathrm{if}\ b_0=0\\
T(x_0^{q-1-p^j}x_1^{b_1+p^j}x_2^{b_2}\cdots
x_n^{b_n}), &\mathrm{if}\ b_0=q-1\\
\big(q(1-T(x_0^{q-1}))-1\big)T(x_1^{b_1+p^j}x_2^{b_2}\cdots
x_n^{b_n}), &\mathrm{if}\ b_0=p^j\\
-J(T^{-p^j},T^{b_0})T(x_0^{b_0-p^j}x_1^{b_1+p^j}x_2^{b_2}\cdots
x_n^{b_n}), &\mathrm{otherwise.}
\end{array}\right.$$
\end{Corollary}
\Proof Let $\rho$ denote the Frobenius automorphism of $R$, which
maps an element of the Teichm\"uller set $T_q$ to its
$p^\mathrm{th}$ power, and let $\mathbf
y=(y_0,\ldots,y_n)=(x_0^{p^j},\ldots,x_n^{p^j})$. We can write
$f_i(\mathbf x)=f_i^{\rho^{-j}}(\mathbf y)=T(y_0^{b_0p^{-j}}
\cdots y_n^{b_np^{-j}})$. We observe that $g_{\ell p^{-j}}$
replaces $y_0$ by $y_0+{\bar\xi}^\ell y_1$. Therefore we can apply
the previous lemma to get
$$g(j)f_i^{\rho^{-j}}(\mathbf y)=\left\{\begin{array}{cl}
0, &\mathrm{if}\ b_0=0\\
T(y_0^{q-2}y_1^{b_1p^{-j}+1}y_2^{b_2p^{-j}}\cdots
y_n^{b_np^{-j}}), &\mathrm{if}\ b_0=q-1\\
\big(q(1-T(y_0^{q-1}))-1\big)T(y_1^{b_1p^{-j}+1}y_2^{b_2p^{-j}}\cdots
y_n^{b_np^{-j}}), &\mathrm{if}\ b_0=p^j\\
-J(T^{-1},T^{b_0p^{-j}})T(y_0^{b_0p^{-j}-1}y_1^{b_1p^{-j}
+1}y_2^{b_2p^{-j}}\cdots y_n^{b_np^{-j}}), &\mathrm{otherwise.}
\end{array}\right.$$
Substituting $\mathbf x$ back in and noting that
$J(\chi^p,\psi^p)=J(\chi,\psi)$ we get the result. \QED

For each basis monomial in $\mathcal M_R$ with at least one
exponent strictly between 0 and $q-1$, we want to construct an
element of $RG$ which acts as the identity on that basis monomial
and annihilates all other members of $\mathcal M_R$. 

\begin{Lemma} \label{torus}
Let $\mathcal M_R=\{ f_1, f_2,\ldots, f_v\}$ be the monomial basis
of $R^{\mathcal L_1}$.  For each $ f_i=T(x_0^{b_0}x_1^{b_1}\cdots
x_n^{b_n}) \in\mathcal M_R$ with some $b_j$ strictly between 0 and
$q-1$, there is an element $h_i\in RG$ with the following
property. If
$$f=c_1 f_1+c_2 f_2+\cdots+c_v f_v$$ is any element in
$R^{\mathcal L_1}$ then
$$h_i f=c_i f_i.$$
\end{Lemma}

\Proof We will construct the required $h_i$ in two steps. 
Let $H$ denote the subgroup of diagonal matrices of $G$.
Then each basis monomial in $\mathcal M_R$ spans a rank one
$RH$-submodule of $R^{\mathcal L_1}$, which is the direct
sum of all such submodules.
Two basis monomials $f_i=T(x_0^{b_0}x_1^{b_1}\cdots x_n^{b_n})$
and $f_j=T(x_0^{b'_0}x_1^{b'_1}\cdots x_n^{b'_n})$
afford the same character of $H$ if and only if
$b_i\equiv b_i'\bmod q-1$ for $0\leq i\leq n$.

Since the order of $H$ is not divisible by $p$, for each
character $\chi$ of $H$, the group ring $RH$ 
contains an idempotent element projecting onto the $\chi$-isotypic
component of $R^{\mathcal L_1}$, the span of all
the basis monomials affording $\chi$.

If none of the exponents the of $f_i$ is divisible by $q-1$,
then no other basis monomials afford the same
character as $f_i$  and we can take $h_i$ to be
the above idempotent. 
Now suppose that some
exponents of $f_i$ are divisible by $q-1$. We proceed successively
for each exponent of $f_i$ which is either 0 or $q-1$. Without
loss of generality, assume that $f_i$ has $b_0=q-1$. We construct
an element $h\in RG$ which annihilates every basis monomial in
${\mathcal M}_R$ for which $b_0=0$ and acts as the identity on $
f_i$. (If $ f_i$ instead has $b_0=0$, then the element we want is
$1-h$.)

Without loss of generality we will take
$b_1=a_{1,t-1}p^{t-1}+\cdots+a_{1,0}$ to be an exponent lying strictly 
between 0 and $q-1$ with $0<a_{1,j}<p-1$ for some $j$. Then we take
the element $h_1=g(j)\in RG$ from Lemma~\ref{frobenius} that
shifts $p^j$ from $b_0$ to $b_1$.  We get
$$h_1f_i=T(x_0^{q-1-p^j}x_1^{b_1+p^j}x_2^{b_2}\cdots x_n^{b^n}).$$
If $e$ is any other basis monomial of the form 
$e=T(x_0^{q-1}x_1^{b_1}\cdots)$, then we similarly have
$$h_1e=T(x_0^{q-1-p^j}x_1^{b_1+p^j}\cdots );$$
and if $x_0$ has exponent 0 in $e$ then from Corollary~\ref{frobenius}, 
we have 
$$h_1e=0.$$
Next we set $h_2=g'(j)\in RG$ to be the analog of $g(j)$ but with
the roles of $x_0$ and $x_1$ interchanged. Noting that here
$b_1+p^j\neq p^j$ (we assumed that $0<b_1<q-1$), we get
\begin{eqnarray*}
h_2h_1f_i&=&-J(T^{-p^j},T^{b_1+p^j})f_i\\
h_2h_1e&=&-J(T^{-p^j},T^{b_1+p^j})e,\quad\textrm{if the exponent
of }x_0 \;{\rm in}\; e\ \mathrm{is}\ q-1\ \\
h_2h_1e&=&0\quad\mathrm{otherwise}.
\end{eqnarray*}
Since there is no carry in the sum $p^j+b_1$, the Jacobi sum 
$J(T^{-p^j},T^{b_1+p^j})$ is a unit in $R$ (cf. Lemma~\ref{jacobidiv}). 
Hence the element $h$ of $RG$ we want is 
$-\frac1{J(T^{-p^j},T^{b_1+p^j})}h_2h_1$.

We can repeat the above process for each exponent of $f_i$ which
is divisible by $q-1$. The product of all the elements we have
constructed is the element $h_i\in RG$ which kills every basis
monomial in ${\mathcal M}_R$ except $f_i$. \QED

We now prove the main result in this section.

\begin{Lemma}\label{basis} Assume $q>2$. There exists an SNF basis of
$R^{\Ll_1}$ for $\eta_{1,r}$, whose reduction 
modulo $\pp$ is $\mathcal M$, and which contains all the basis monomials of ${\mathcal M}_R$ having at least one exponent lying strictly between $0$ and $q-1$.
\end{Lemma}

\Proof By Corollary~\ref{snfbasis}, there exists an SNF basis 
$\bb=\cup_{j=0}^{\ell -1}\bb_j$ of 
$R^{\Ll_1}$ for $\eta_{1,r}$ such that the reduction of $\bb$ modulo 
$\pp$ is $\mathcal M$.  
Let $f\in\bb$, and let the reduction of $f$ modulo $\pp$ be 
$$\overline{f}=x_0^{b_0}x_1^{b_1}\cdots x_n^{b_n}\in \mathcal M,$$
with some $b_j$ satisfying $0<b_j<q-1$.
Let ${\mathcal M}_R=\{f_1, f_2, \ldots ,f_v\}$ with $f_1=T(x_0^{b_0}x_1^{b_1}\cdots x_n^{b_n})$, where $v=|\Ll_1|$. We write 
$$f=c_1f_1+c_2f_2+\cdots +c_vf_v, \; c_i\in R.$$
Since $\overline{f}=\overline{f_1}$, we see that $\overline{c_1}=1$, hence 
$c_1$ is a unit in $R$. Since there is an exponent $b_j$ lying 
strictly between $0$ and $q-1$, by Lemma~\ref{torus}, 
we can find $h_1\in RG$ such that $h_1f=c_1f_1$. 
In the notation of Definition~\ref{SNFbasis} with $M=R^{\Ll_1}$, 
we see that if $f\in\bb_j$, then $f_1\in M_j$ 
since $M_j$ is an $RG$-submodule , so 
$\bb'=(\bb\setminus\{f\})\cup\{f_1\}$ 
is again an SNF basis of $R^{\Ll_1}$ for $\eta_{1,r}$. We can repeat this process for every element in $\bb$ whose reduction modulo $\pp$ has one exponent strictly lying between $0$ and $q-1$. At the end, we obtain the required SNF basis of $R^{\Ll_1}$ for $\eta_{1,r}$.
\QED

We will use ${\mathcal M}_R'$ to denote the special SNF basis of $R^{\Ll_1}$ for $\eta_{1,r}$ produced by Lemma~\ref{basis}. Again 
{\it the type of $f\in {\mathcal M}_R'$} is defined to be that of 
$\overline{f}\in {\mathcal M}$.

\begin{Lemma}\label{cyclic} 
The invariants of  $\eta_{1,r}$
corresponding to two elements of ${\mathcal M}_R'$ of types
$(s_0,\ldots,s_{t-1})$ and $(s_1,s_2,\ldots,s_{t-1},s_0)$,
respectively, are equal.
\end{Lemma}
\Proof 
We may assume $t\geq2$ since there is nothing to prove
otherwise. 
For any type $\xi\in {\mathcal H}$, we can always find a basis 
monomial $f\in {\mathcal M}_R$ of type $\xi$ and with at least one exponent 
lying strictly between 0 and $q-1$. Hence $f\in {\mathcal M}_R'$. 
By Corollary~\ref{types}, the invariants of $\eta_{1,r}$ corresponding to 
two elements in ${\mathcal M}_R'$ of the same type are equal. Therefore we 
may assume that the two elements of ${\mathcal M}_R'$ in the statement of 
the lemma are actually in ${\mathcal M}_R$. 

The Frobenius field automorphism
\begin{eqnarray*}
\rho:\ \Ff_q &\rightarrow & \Ff_q\\
x_i &\mapsto & x_i^p\\
\end{eqnarray*}
applied to the coordinates of $V$ is an automorphism of the
projective geometry. It maps points to points, subspaces to
subspaces, and preserves incidence . The image of a point
$Z=(x_0,\ldots,x_n)$ is $Z^\rho=(x_0^p,\ldots,x_n^p)$
and for an $r$-subspace $Y$,
$Y^\rho$ is the $r$-subspace containing the images of all the
points incident with $Y$. Given a monomial function
$f_i=T(x_0^{b_0}\cdots x_n^{b_n})$ we have
$$f_i^\rho=T(x_0^{pb_0}\cdots x_n^{pb_n}).$$

Clearly if $f_i$ is
of type $(s_0,\ldots, s_{t-1})$ then $f_i^\rho$ is of type
$(s_{t-1},s_0,\ldots,s_{t-2})$ because $\lambda_j$ becomes
$\lambda_{j+1}$ in (\ref{lambda}). It is also clear that
$$f_i(Z^\rho)=f_i^\rho(Z)$$
so that
$$\eta_{1,r}(f_i)(Y^\rho)=\eta_{1,r}(f_i^\rho)(Y).$$
As $Y$ runs through $R^{\mathcal L_r}$ so does $Y^\rho$. 
Thus, the coordinates of $\eta_{1,r}(f_i)$ are the
same as the coordinates of $\eta_{1,r}({f_i}^\rho)$ but permuted
by $\rho$, so the invariants corresponding to
$f_i$ and $f_i^\rho$ are equal.
\QED

\section{The Proof of Theorem~\ref{main}}

Our aim in this section to prove Theorem~\ref{main} and we will 
achieve this by proving the more detailed result Theorem~\ref{pvalue} below. 
Our proof depends on Lemma~\ref{lower}, which gives lower bounds on 
the $p$-adic valuations of the coordinates of $\eta_{1,r}(f)$, where 
$f\in {\mathcal M}_R$, and the results in Section 5 and 6. 

We first prove a lemma.

\begin{Lemma}\label{prank}
Let $f$ be a nonconstant basis monomial in ${\mathcal M}_R$.
Then $p$ does not divide $\eta_{1,r}(f)$ if and
only if  $f$ has type $(s_0,s_1,\ldots,s_{t-1})$, 
with $s_j\ge r$ for all $0\le j\le t-1$.
\end{Lemma}

\Proof 
Let $\overline f$ be the image modulo $\pp$ of $f$.
Then $p$ does not divide $\eta_{1,r}(f)$ if and only if the image
of $\overline f$ under the induced map 
$\overline\eta_{1,r}:\Ff_q^{\Ll_1}\rightarrow\Ff_q^{\Ll_r}$
is nonzero. 
Suppose that $s_j<r$ for some $j$. 
By Lemma~\ref{lower}, $p|\eta_{1,r}(f)$.
That is, only those basis monomials
$\overline f$ of type $(s_0,s_1,\ldots ,s_{t-1})$, where $s_j\ge
r$ for all $0\le j\le t-1$, could possibly have nonzero image under 
$\overline\eta_{1,r}$. On the other hand, by Hamada's formula, 
rank of $\overline\eta_{1,r}$ is
equal to one plus the number of $\overline f$'s with this property.
Therefore, the images of all such basis monomials
must be linearly independent, in particular, nonzero.
Hence $p\not{|}\eta_{1,r}(f)$, if and only if $f$ has type
$(s_0,s_1,\ldots,s_{t-1})$, where $s_j\ge r$ for all $0\le j\le t-1$.
This completes the proof.\QED

\begin{Theorem} Let $\mathcal M'_R=\{f'_1,f'_2,\ldots ,f'_v\}$ with $\bar f'_1=1\in \mathcal M$. Let the type of $f'_i$, $2\leq i\leq v$, be $(s_0^{(i)},s_1^{(i)},\ldots ,s_{t-1}^{(i)})$ and let $p^{\beta_i}$ be the invariant of $\eta_{1,r}$ corresponding to $f'_i$. Then
$$
\beta_i=\sum_{j=0}^{t-1}\max\{0,r-s_j^{(i)}\}.
$$
\label{pvalue}
\end{Theorem}

\Proof We shall assume that $t\geq 2$. When $t=1$ a similar and easier 
argument works, but we omit the details to keep the notation simple and 
the argument clear, since this case is already known \cite{sinp}.
Let $\alpha_i=\sum_{j=0}^{t-1}\max\{0,r-s_j^{(i)}\}$ and let 
$f_i\in\mathcal M_R$ be the basis monomial which
has the same reduction modulo $\pp$ as $f'_i$, namely
$f_i=T(\bar f_i')$. 
We use the notation of Definition~\ref{SNFbasis} with $M=R^{\Ll_1}$
and $\phi=\eta_{1,r}$.
By Lemma~\ref{lower}, we have $f_i\in M_{\alpha_i}$.
Since the image of $\bar f_i=\bar f'_i$
in $\overline M_{\beta_i}/\overline M_{\beta_i+1}$ is not zero,
it follows that $\alpha_i\leq \beta_i$.
 
Suppose by way of contradiction that $\beta_k>\alpha_k$ for some $k$. 
Let $f_k=T(x_0^{b_0}x_1^{b_1}\cdots x_n^{b_n})$ be of type $(s_0,
s_1,\ldots,s_{t-1})$ (here we suppressed the superscript $(k)$ of $s_j$ 
to keep the notation simple). Assume that we have picked $k$ so that if
$\alpha_j<\alpha_k$ then $\alpha_j=\beta_j$. By Lemma~\ref{cyclic}
we can assume for convenience that
$s_1=\min\{s_0,\ldots,s_{t-1}\}$. We have
$$\lambda_0=ps_1-s_0\le n(p-1)$$
with equality only if $s_0=s_1=\cdots=s_{t-1}=n$ and
$$\lambda_1=ps_2-s_1\ge1.$$

We note that the case $s_0=s_1=\cdots=s_{t-1}=n$ will not occur by
our assumption that $\beta_k>\alpha_k$. The reason is as follows.
If $f_k$ has type $(s_0,s_1,\ldots ,s_{t-1})=(n,n,\ldots ,n)$, by
Lemma~\ref{prank}, we see that $p\not\hspace{-0.03in}{|}\eta_{1,r}(f_k)$. Since 
$\bar f_k'=\bar f_k$, we have $p\not\hspace{-0.02in}{|}\eta_{1,r}(f'_k)$. But the
invariant corresponding to $f'_k$ is $p^{\beta_k}$, and we assumed
that $\beta_k>\alpha_k=0$, so $p|\eta_{1,r}(f'_k)$, a contradiction.

By Lemma~\ref{types}, basis vectors in ${\mathcal M}'_R$ 
of the same type correspond to the same invariant, so in the sum
$\lambda_0=\sum_{i=0}^na_{i,0}$ we can assume that $a_{0,0}=0$, and 
we can also assume that $a _{1,0}<p-1$ since the case  
$s_0=s_1=\cdots=s_{t-1}=n$ has been excluded. In the sum 
$\lambda_1=\sum_{i=0}^na_{i,1}$, we can assume that $a_{0,1}\ge1$. 
By these assumptions, we see that $0<p\leq b_0<q-1$, hence from our 
definition of ${\mathcal M}'_R$ we have
$$f'_k=f_k.$$

Since the exponent $b_0$ in $f_k$ is not equal to 1, applying the
group ring element $h\in RG$ in Lemma~\ref{specialg}, we get
\begin{equation}\label{action}
hf'_k=hf_k=-J(T^{-1},T^{b_0})T(x_0^{b_0-1}x_1^{b_1+1}x_2^{b_2}\cdots
x_n^{b_n}).
\end{equation}
Set $T(x_0^{b_0-1}x_1^{b_1+1}x_2^{b_2}\cdots
x_n^{b_n}):=f_{\ell}\in {\mathcal M}_R$. The type of $f_{\ell}$ is
$(s_0,s_1+1,s_2,\ldots,s_{t-1})$ because we have increased
$\lambda_0$ by $p$ and decreased $\lambda_1$ by 1. Also note that
$b_0-1$ is still strictly between 0 and $q-1$, so
$f_{\ell}=f'_{\ell}\in {\mathcal M}'_R$. As for the coefficient of
$f_{\ell}$ in (\ref{action}), Lemma~\ref{jacobidiv} tells us
that $p$ divides $J(T^{-1},T^{b_0})$ exactly once because when 1
is added to $q-1-b_0$ there is exactly one carry: from the ones
place to the $p$ place of the sum. Since
$p^{\beta_k}|\eta_{1,r}(f'_k)$ and $\eta_{1,r}$ is an
$RG$-module homomorphism, we have
$$p^{\beta_k}|\eta_{1,r}(hf'_k).$$
Since $p\parallel J(T^{-1},T^{b_0})$, we get
$$p^{(\beta_k-1)}\mid\eta_{1,r}(f'_{\ell}),$$
where the type of $f'_{\ell}$ is $(s_0,s_1+1,s_2,\ldots,s_{t-1})$.
Since we assumed that $\alpha_k$ was the smallest such that
$\alpha_k<\beta_k$, we must conclude that
$$\sum_{j=0}^{t-1}\max\{0,r-s_j\}=\sum_{j=0, j\neq
1}^{t-1}\max\{0,r-s_j\}+\max\{0, r-(s_1+1)\}.$$
That is, $s_1\geq r$. Since $s_1$ is assumed to be the smallest among
$s_j, 0\leq j\leq t-1$, we see that
$$s_j\ge r, \;  0\le j\le t-1, \; {\rm and}\; {\rm hence}\;  \alpha_k=0.$$
By Lemma~\ref{prank}, $p\not{|}\eta_{1,r}(f_k)$, so $p\not{|}\eta_{1,r}(f'_k)$
since $f'_k=f_k$. However we have assumed that $\beta_k>\alpha_k=0$, that is,
$p|\eta_{1,r}(f'_k)$. This is a contradiction. The theorem is
proved.\QED

The following corollary is immediate.

%theorem shows that the bound in Theorem~\ref{wanlow} is exact.

\begin{Corollary}\label{snfmr} The monomial basis $\mathcal M_R$ is
an SNF basis of $R^{\Ll_1}$ for the map $\eta_{1,r}$ and
the invariant of $\eta_{1,r}$ corresponding to a monomial of type 
$(s_0,\ldots,s_{t-1})$ is equal to
$$
\sum_{j=0}^{t-1}\max\{0,r-s_j\}.
$$

\end{Corollary}

\begin{Remark}{\rm  We have seen that, for each $r$,
the $R\mathrm{GL}(n+1,q)$ homomorphism $\eta_{1,r}$
defines a filtration $\{\overline M_i\}$ 
of $\Ff_q^{\Ll_1}$ by 
$\Ff_q\mathrm{GL}(n+1,q)$-modules. In the case $r=n$, it follows from
Theorem~\ref{pvalue} and \cite{bsin}, Theorems A and B,
that this filtration is equal to the radical filtration,
the most rapidly descending filtration with semisimple
factors. Equivalently, $M_i=J^i(\Ff_q^{\Ll_1})$, where
$J$ is the Jacobson radical of the group algebra $\Ff_q\mathrm{GL}(n+1,q)$.}
\end{Remark}

\section{The Invariant Factors of the Incidence between points and $r$-flats in $\mathrm{AG}(n,q)$}

In this section, we consider the incidence between points and
$r$-flats in the affine geometry $\mathrm{AG}(n,q)$. We will view
$\mathrm{AG}(n,q)$ as obtained from $\mathrm{PG}(n,q)$ by deleting
a hyperplane and all the subspaces it contains. Let $H_0$ be the
hyperplane of $\mathrm{PG}(n,q)$ given by the equation $x_0=0$.
Then for any integer $r$, $0\leq r\leq n$, the set of $r$-flats of
$\mathrm{AG}(n,q)$ is
$$\mathcal F_r=\{Y\setminus (Y\cap H_0)\mid Y\in \Ll_{r+1}\}.$$
(The empty set is not considered as an $r$-flat for any $r$.) In
particular, the set of points of $\mathrm{AG}(n,q)$ is $\mathcal
F_0$. We define the incidence map
\begin{equation}\label{defetaaff}
\eta'_{0,r}: \Zz^{\mathcal F_0}\rightarrow \Zz^{\mathcal F_r}
\end{equation}
by letting $\eta'_{0,r}(Z)=\sum_{Y\in\mathcal F_r, Z\subset Y}Y$
for every $Z\in\mathcal F_0$, and then extending $\eta'_{0,r}$
linearly to $\Zz^{\mathcal F_0}$. Similarly, we define
$\eta'_{r,0}$ to be the map from $\Zz^{\mathcal F_r}$ to
$\Zz^{\mathcal F_0}$ sending an $r$-flat of ${\rm AG}(n,q)$ to the
formal sum of all points incident with it. Let $A_1$ be the matrix
of $\eta'_{0,r}$ with respect to the standard bases of
$\Zz^{\mathcal F_0}$ and $\Zz^{\mathcal F_r}$. We have the
following counterpart of Theorem~\ref{p'part}.

\begin{Theorem}
The invariant factors of $A_1$ are all powers of $p$.
\end{Theorem}

\Proof The proof is parallel to that of Theorem~\ref{p'part}. We
will actually work with $A_1^{\top}$, which is the matrix of
$\eta'_{r,0}:\Zz^{\mathcal F_r}\rightarrow \Zz^{\mathcal F_0}$
with respect to the standard bases of $\Zz^{\mathcal F_r}$ and
$\Zz^{\mathcal F_0}$. 
We define $$\epsilon':\Zz^{\mathcal F_0}\rightarrow\Zz$$ to be the function sending each element in $\Fl_0$ to 1. Clearly $\epsilon'$ maps
$\Zz^{\mathcal F_0}$ onto $\Zz$ and $\im\eta'_{r,0}$
onto $q^r\Zz$. Thus,
$\Zz^{\mathcal F_0}/(\Ker\epsilon'+\im\eta'_{r,0})\cong\Zz/q^r\Zz$,
and we are reduced to proving that
$(\Ker\epsilon'+\im\eta'_{r,0})/\im\eta'_{r,0}$
is a $p$-group. The proof goes in exactly the same way as
that in Theorem~\ref{p'part}. Note that $\mathrm{Ker}(\epsilon')$
is spanned by elements in $\Zz^{\Fl_0}$ of the form $u-w$, where
$u$ and $w$ are distinct points of ${\rm AG}(n,q)$; so it is
enough to show that $q^r(u-w)\in\im(\eta'_{r,0})$ for any
two distinct points $u$ and $w$. We pick an $(r+1)$-flat
containing the two distinct points $u$ and $w$ and let
$\tilde\eta'_{0,r}$ be the restricted map. The number of $r$-flats
through one point in $\AG(r+1,q)$ is $(q^{r+1}-1)/(q-1)$
while the number of $r$-flats through two points in
$\AG(r+1,q)$ is $(q^{r}-1)/(q-1)$ so we get
$$\eta'_{r,0}(\tilde\eta'_{0,r}(z))=q^{r}z+
\frac{q^{r}-1}{q-1}{\mathbf j}_U$$ for any point $z$. Therefore
$$\eta'_{r,0}(\tilde\eta'_{0,r}(u-w))=q^r(u-w).$$
This completes the proof. \QED

In view of the above theorem, we view $A_1$ as a matrix with
entries from $R=\Zz_p[\xi_{q-1}]$. The Smith normal form of $A_1$
over $R$ will completely determine the Smith normal form of $A_1$
over $\Zz$. We will get the $p$-adic invariants of $A_1$ from the
invariants of the incidence between points and $(r+1)$-spaces in
${\rm PG}(n,q)$ and those of the incidence between points and
$(r+1)$-spaces in ${\rm PG}(n-1,q)$.

Let $A$ be the matrix of the incidence map
$\eta_{1,r+1}:R^{\Ll_1}\rightarrow R^{\Ll_{r+1}}$ with respect to
the standard bases of $R^{\Ll_1}$ and $R^{\Ll_{r+1}}$. We want to
partition $A$ into a certain block form. For this purpose, we
define
$$\Ll_{1}^{H_0}=\{Z\in \Ll_{1}\mid Z\subseteq H_0\},$$
and
$$\Ll_{r+1}^{H_0}=\{Y\in \Ll_{r+1}\mid Y\subseteq H_0\}.$$
So we have the partitions
$$\Ll_{1}=\Fl_0\cup \Ll_1^{H_0},$$
and
$$\Ll_{r+1}=\Fl_r\cup \Ll_{r+1}^{H_0}.$$
We now partition $A$ as

\[
A = 
\left[
\begin{array}{c|c}
\makebox(0,0)[bl]{\hspace*{-0.75ex}$\overbrace{\phantom{A_1}}^{\Fl_0}$} A_1 &
\makebox(0,0)[bl]{\hspace*{-0.75ex}$\overbrace{\phantom{A_2}}^{\Ll_1^{H_0}}$} A_2\\
\hline
0 & A_3
\end{array}
\right]\hspace*{-4ex}
\begin{array}{c}
\left.\phantom{A_1}\right\}\makebox(0,0)[bl]{\scriptsize$\Fl_r$}\\
%\phantom{\hline}
\left.\phantom{A_1}\right\}\makebox(0,0)[cl]{\scriptsize$\Ll_{r+1}^{H_0}$}
\end{array}
\]
where $A_3$ is the incidence matrix of the incidence between
$\Ll_1^{H_0}$ and $\Ll_{r+1}^{H_0}$, which can be thought as the
matrix of the incidence between points and $(r+1)$-spaces in
${\PG}(n-1,q)$.

In order to obtain the SNF of $A_1$, we need to modify the
monomial basis ${\mathcal M}_R$ of $R^{\Ll_1}$ slightly. We
replace the constant monomial in ${\mathcal M}_R$ by
$T(x_0^{q-1}x_1^{q-1}\cdots x_n^{q-1})$ and denote the resulting
set by ${\mathcal M}^*_R$. Note that ${\mathcal M}^*_R$ is still a
basis of $R^{\Ll_1}$ because
$(1-a_0^{q-1})(1-a_1^{q-1})\cdots(1-a_n^{q-1})=0$ for each point
$(a_0,a_1,\ldots ,a_n)$ of $\mathrm{PG}(n,q)$. Furthermore 
${\mathcal M}^*_R$ is an SNF basis of $R^{\Ll_1}$ for $\eta_{1,r+1}$ since 
${\mathcal M}_R$ is an SNF basis of $R^{\Ll_1}$ for $\eta_{1,r+1}$ and 
the invariant corresponding to $T(x_0^{q-1}x_1^{q-1}\cdots x_n^{q-1})$ is 1. 
So we have the factorization
\begin{equation}\label{smithequ}
P^*D=AQ^*, \end{equation}
where the columns of $Q^*$ are the basis
vectors in $M_R^*$ written with respect to the standard basis of
$R^{\Ll_{1}}$, $P^*$ is nonsingular over $R$ and $D$ is the Smith
normal form of $A$.

We now partition ${\mathcal M}^*_R$ as $\bb_1\cup\bb_2$, where
$$\bb_1=\{T(x_0^{b_0}x_1^{b_1}\cdots x_n^{b_n})\mid b_0\neq 0, T(x_0^{b_0}x_1^{b_1}\cdots
x_n^{b_n})\in {\mathcal M}^*_R\},$$ and
$$\bb_2=\{T(x_0^{b_0}x_1^{b_1}\cdots x_n^{b_n})\mid b_0=0, T(x_0^{b_0}x_1^{b_1}\cdots
x_n^{b_n})\in {\mathcal M}^*_R\}.$$ We partition the matrix $Q^*$
according to the partition of ${\mathcal M}^*_R$ as
$\bb_1\cup\bb_2$ and the partition of $\Ll_1$ as
$\Fl_0\cup\Ll_1^{H_0}$. Explicitly we have

\[
Q^* = 
\left[
\begin{array}{c|c}
\makebox(0,0)[bl]{\hspace*{-0.75ex}$\overbrace{\phantom{Q_1}}^{\bb_1}$} Q_1 &
\makebox(0,0)[bl]{\hspace*{-0.75ex}$\overbrace{\phantom{Q_2}}^{\bb_2}$} Q_2\\
\hline
0 & Q_3
\end{array}
\right]\hspace*{-4ex}
\begin{array}{c}
\left.\phantom{A_1}\right\}\makebox(0,0)[bl]{\scriptsize$\Fl_0$}\\
%\phantom{\hline}
\left.\phantom{A_1}\right\}\makebox(0,0)[cl]{\scriptsize$\Ll_1^{H_0}$}
\end{array}
\]
where the 
columns of $Q_3$ are the basis vectors in $\{f|_{\scriptscriptstyle H_0}\mid
f\in\bb_2\}$ written with respect to the standard basis of
$R^{\Ll_1^{H_0}}$.

Now we rewrite (\ref{smithequ}) according to the block forms of
the matrices $A$ and $Q^*$. We have
\begin{equation}\label{partition}
\left(\begin{array}{ccc}P_1&P_3&P_5\\0&P_2&P_4\end{array}\right)
\left(\begin{array}{cc}D_1&0\\0&D_2\\0&0\end{array}\right)=
\left(\begin{array}{cc}A_1&A_2\\0&A_3\end{array}\right)
\left(\begin{array}{cc}Q_1&Q_2\\0&Q_3\end{array}\right)
\end{equation}
which gives us
$$P_1D_1=A_1Q_1,$$
and
$$P_2D_2=A_3Q_3.$$
Since $P_1$ and $Q_1$ inherit the property that the reductions
modulo $\pp$ of their columns are linearly independent, $D_1$ must
be the Smith normal form of $A_1$. By Corollary~\ref{snfmr}, 
$\{f|_{\scriptscriptstyle H_0}\mid f\in\bb_2\}$ is an SNF basis of $R^{\Ll_1^{\scriptscriptstyle H_0}}$
for the incidence map $\eta_{1,r+1}$ between points and $(r+1)$-spaces in ${\rm PG}(n-1,q)$. We see that $D_2$ is the Smith normal form of $A_3$.

For any $n\geq 2$, $1<i\leq n$, and $\alpha\geq 0$, let
$m(\alpha,n,i)$ denote the multiplicity of $p^\alpha$ as a
$p$-adic invariant of the incidence between points and projective
$(i-1)$-dimensional subspaces in $\mathrm{PG}(n,q)$. (The numbers
$m(\alpha,n,i)$ are determined by Theorem~\ref{main}.) We have the
following theorem.

\begin{Theorem} The $p$-adic invariants of $A_1$ are $p^{\alpha}$,
$0\leq\alpha\leq rt$, with multiplicity
$m(\alpha,n,r+1)-m(\alpha,n-1,r+1)$.
\end{Theorem}

\Proof From (\ref{partition}), we see that the multiplicity of
$p^{\alpha}$ as invariant of $A_1$ is equal to the number of times
$p^\alpha$ appearing in $D$ minus the number of times $p^{\alpha}$
appearing in $D_2$. \QED

\end{document}